\documentclass[letterpaper, 10 pt, conference]{ieeeconf}  
\IEEEoverridecommandlockouts                              

\usepackage{graphicx}
\usepackage{amssymb,mathtools,amsfonts}
\usepackage{hyperref}
\usepackage{multirow}
\usepackage{color}
\usepackage{psfrag}
\usepackage{empheq}
\usepackage{tabularx}
\usepackage{bm} 
\usepackage{url}
\usepackage{float}
\usepackage{comment}  
\usepackage{caption, subcaption}
\usepackage{alphalph}
  
\usepackage{amsthm}

\newcommand{\R}{{\mathbb{R}}}
\newcommand{\dd}{{\mathrm{d}}}

\newcommand{\Nsys}{{n_{f}}} 
\newcommand{\Nstg}{{n_\mathrm{s}}} 
\newcommand{\NFE}{N_\mathrm{FE}} 

\newcommand{\C}{{\mathcal{C}}}
\newcommand{\J}{{\mathcal{J}}}

\newcommand{\lambdan}{\lambda^{\mathrm{n}}}
\newcommand{\lambdap}{\lambda^{\mathrm{p}}}

\newcommand{\Lambdan}{\Lambda^{\mathrm{n}}}
\newcommand{\Lambdap}{\Lambda^{\mathrm{p}}}

\newcommand{\irk}{{\mathrm{rk}}}

\newcommand{\fesd}{{\mathrm{fesd}}}

\newcommand{\Nctrl}{N}

\newcommand{\ts}{{t_{\mathrm{s}}}}

\newtheorem{theorem}{Theorem}
\newtheorem{example}{Example}
\newtheorem{remark}[theorem]{Remark}

\newtheorem{proposition}[theorem]{Proposition}

\newtheorem{lemma}[theorem]{Lemma}

\begin{document}

\title{Finite Elements with Switch Detection for Direct Optimal Control of Nonsmooth Systems with Set-Valued Step Functions}
\author{Armin Nurkanovi\'c$^{1}$, Jonathan Frey$^{1,2}$, Anton Pozharskiy$^{1}$, Moritz Diehl$^{1,2}$
\thanks{This research was supported by the DFG via Research Unit FOR 2401 and project 424107692 and by the EU via ELO-X 953348.}
\thanks{$^1$Department of Microsystems Engineering (IMTEK), 	
$^2$Department of Mathematics, University of Freiburg, Germany,  
\sloppy\texttt{\{armin.nurkanovic,jonathan.frey,moritz.diehl\} @imtek.uni-freiburg.de, anton.pozharskiy@ merkur.uni-freiburg.de}
		}
}
\maketitle
\thispagestyle{empty} 
\begin{abstract}
This paper extends the Finite Elements with Switch Detection (FESD) method~\cite{Nurkanovic2022} to optimal control problems with nonsmooth systems involving set-valued step functions.
Logical relations and common nonsmooth functions within a dynamical system can be expressed using linear and nonlinear expressions involving step functions.
A prominent subclass of these systems are Filippov systems.
The set-valued step function can be expressed by the solution map of a linear program, and using its KKT conditions allows one to transform the initial system into an equivalent dynamic complementarity system (DCS).
Standard Runge-Kutta (RK) methods applied to DCS have only first-order accuracy.
The FESD discretization makes the step sizes degrees of freedom and adds further constraints that ensure exact switch detection to recover the high-accuracy properties that RK methods have for smooth ODEs.
We use the novel FESD method for the direct transcription of optimal control problems.
All methods and examples in this paper are implemented in the open-source software package~\texttt{NOSNOC}.
\end{abstract}
\section{Introduction}
In this paper, we introduce a high-accuracy method for discretizing and solving nonsmooth Optimal Control Problems (OCPs) of the following form:
\begin{subequations} \label{eq:ocp}
	\begin{align}
		\min_{x(\cdot),u(\cdot)} \quad & \int_{0}^{T} L(x(t),u(t))\dd t +  R(x(T)) \\
		\textrm{s.t.} \quad  x_{0} &= s_0, \\
		\dot{x}(t) &\! \in F(x(t), u(t), \Gamma(c(x(t)))),\ \textrm{a.a. } t \in [0, T], \label{eq:ocp_di} \\
		0&\geq G_{\mathrm{p}}(x(t),u(t)),\ t \in [0,T],\\
		0&\geq G_{\mathrm{t}}(x(T)),
	\end{align}
\end{subequations}
where $L: \R^{n_x} \times \R^{n_u} \to \R$ is the running cost and $R:\R^{n_x}\to \R$ is the terminal cost, $s_0\in\R^{n_x}$ is a given initial value.
The right-hand side of the Differential Inclusion (DI) in Eq. \eqref{eq:ocp_di} is the set-valued mapping 
$F: \R^{n_x} \times \R^{n_u} \times \R^{n_c} \to \mathcal{P}(\R^{n_x})$.
The path and terminal constraints are defined by the functions $G_{\mathrm{p}} : \R^{n_x}  \times \R^{n_u} \to \R^{n_{\mathrm{p}}}$ and $G_{\mathrm{t}} : \R^{n_x}  \to \R^{n_{\mathrm{t}}}$, respectively.

The OCP is nonsmooth due to the DI in Eq. \eqref{eq:ocp_di}. 
The function $c(x) \in \R^{n_c}$ contains $n_c$ \textit{switching functions}.
The set-valued function $\Gamma: \R^{n_c} \to \mathcal{P}(\R^{n_c})$ is defined as the concatenation of scalar step functions, i.e., for $y \in \R^{n_c}$ we have
$\Gamma(y)= [\gamma(y_1),\dots,\gamma(y_{n_c})]^\top \in \R^{n_c}$, where $\gamma : \R \to \mathcal{P}(\R)$ is defined as:
\begin{align}\label{eq:step_function}
	\gamma(y_i) &= \begin{cases}
		\{1\},\ &y_i>0, \\
		[0,1],\ &y_i=0, \\
		\{0\},\ &y_i<0.
	\end{cases}
\end{align}
Note that there are no particular restrictions on how the components $\Gamma(c(x))$ enter the r.h.s. of the DI \eqref{eq:ocp_di}.
This DI is an instance of so-called Aizerman–Pyatnitskii DIs., cf. \cite[page 55, Definition c]{Filippov1988}.
A prominent and well-studied subclass of these DIs, on which we focus in the sequel, are Filippov DIs. 

Set-valued step functions provide an intuitive way to model logical \textit{if-else} and \textit{and-or} relations in dynamical systems.
In addition, several other common nonsmooth functions, such as $\mathrm{sign},\min$, and $\max$, can be easily expressed via step functions.
Smoothed versions of the step function are often used in the numerical simulation of Piecewise Smooth Systems (PSS)~\cite{Guglielmi2022,Machina2011}.
Step functions are often used in the modeling of gene regulatory networks~\cite{Acary2014,Machina2011}.
Another application is to express Filippov sets in sliding modes on surfaces with co-dimension higher than one~\cite{Dieci2011}.
Moreover, several classes of systems with state jumps can be reformulated into PSS via the time-freezing reformulation~\cite{Nurkanovic2023a,Nurkanovic2021,Nurkanovic2022a,Halm2021}.
Thus, the formulation \eqref{eq:ocp} covers a wide class of practical problems.

The OCP \eqref{eq:ocp} is difficult to solve numerically for several reasons.
Direct methods first discretize the infinite-dimensional OCP~\eqref{eq:ocp} and solve then a finite-dimensional Nonlinear Program (NLP).
However, the discretization of a DI \eqref{eq:ocp_di} presents several pitfalls.
First, standard time-stepping methods for DIs have only first-order accuracy~\cite{Acary2010}.
Stewart and Anitescu~\cite{Stewart2010} have shown that the numerical sensitivities obtained from standard time-stepping are incorrect regardless of the integrator step size.
Moreover, the sensitivities of smoothed approximations of \eqref{eq:ocp_di} are correct only under the very restrictive assumption that the step size is sufficiently smaller than the smoothing parameter. 
Smoothing and wrong sensitivities can lead to artificial local minima and jeopardize the progress of NLP solvers~\cite{Nurkanovic2020}.
Furthermore, the discretized OCPs are nonsmooth NLPs. 
In summary, even for moderate accuracy, many optimization variables and a huge computational load are required.

Some of these drawbacks are overcome by the recently introduced Finite Elements with Switch Detection (FESD) method~\cite{Nurkanovic2022}. 
FESD was originally developed for Filippov DIs, which are transformed into equivalent Dynamic Complementarity System (DCS) via Stewart's reformulation~\cite{Stewart1990a}.
This method starts with a standard Runge-Kutta (RK) discretization of the DCS and, inspired by \cite{Baumrucker2009}, allows the integrator step sizes to be degrees of freedom. 
Additional constraints ensure exact switch detection (and thus higher-order accuracy) and correct computation of numerical sensitivities (and avoid convergence to spurious solutions). 
This overcomes the aforementioned fundamental limitations of standard time-stepping discretization methods.

\paragraph*{Contributions}
In this paper, we extend the FESD method to DIs (and the associated OCP \eqref{eq:ocp}) governed by set-valued step functions.
Our approach is to express the step function as the solution map of a linear program, and using its KKT conditions, we write the DI into an equivalent DCS. 
The DCS is discretized with a standard RK method.
However, we let the step sizes be degrees of freedom and introduce additional equations that ensure exact switch detection.
This recovers the high-accuracy properties that RK methods have for smooth ODEs.
FESD results in Mathematical Programs with Complementarity Constraints (MPCC), which can be efficiently solved in a homotopy loop with off-the-shelf NLP solvers~\cite{Scholtes2001,Anitescu2007}.
We illustrate the efficacy of the new formulation on numerical simulation and OCP examples.
All methods and examples of this paper are implemented in the open-source package~\texttt{NOSNOC}~\cite{Nurkanovic2022b}\footnote{MATLAB: \url{https://github.com/nurkanovic/nosnoc},\\Python: \url{https://github.com/FreyJo/nosnoc_py}}.
\paragraph*{Notation}
The complementary conditions for two vectors  $x,y \in \R^{n}$ {read as ${0\leq x\perp y\geq 0}$, where $x \perp y$ means $x^{\top}y =0$}.
For two scalar variables $a,b$ the so-called C-functions have the property $\phi(a,b) = 0 \iff a\geq 0, b\geq 0, ab = 0$.
A famous example is the Fischer-Burmeister function $\phi_{\mathrm{FB}}(a,b) = a+b-\sqrt{a^2+b^2}$.
If $x,y \in \R^{n}$, we use $\phi(\cdot)$ component-wise and define $\Phi(x,y) = (\phi(x_1,y_1),\dots,\phi(x_{n},y_{n}))$.
The concatenation of two column vectors $x\in \R^{n}$, $y\in \R^{m}$ is denoted by $(x,y)\coloneqq[x^\top,y^\top]^\top$. 
Given a matrix $S \in \R^{n \times m}$, its $i$-th row is denoted by $S_{i,\bullet}$ and its $j$-th column is denoted by $S_{\bullet,j}$.

\section{Filippov systems and the equivalent dynamic complementarity system}\label{sec:pss_and_dcs}
In this section, we show how the Filippov convexification of a PSS can be expressed via step functions, and how to transform this system into an equivalent DCS.
\subsection{Filippov systems via step functions}\label{sec:step_dcs}
We focus on the PSS systems and their Filippov convexification as the most prominent representative of DIs with set-valued step functions.
Most developments are straightforwardly generalized.
A controlled PSS is defined as
\begin{align} \label{eq:pss}
	\dot{x}(t) = f_i(x(t),u(t)),\ &\mathrm{if} \  x(t) \in R_i \subset \R^{n_x}, i \in \mathcal{J},
\end{align}
where $R_i$ are disjoint, connected, and open sets, and $\mathcal{J} \coloneqq \{ 1,\dots,\Nsys  \}$.
The sets $R_i$ are assumed to be nonempty and to have piecewise-smooth boundaries $\partial R_i$.
It is assumed that $\overline{\bigcup\limits_{i\in \mathcal{J}} R_i} = \R^{n_x}$, and that $\R^{n_x} \setminus \bigcup\limits_{i\in\mathcal{J}} R_i$ is a set of measure zero.
The functions $f_i(\cdot)$ are assumed to be Lipschitz and at least twice continuously differentiable functions on an open neighborhood of $\overline{R}_i$.
Here, $u(t)$ is a sufficiently regular externally chosen control function, e.g., obtained as a solution to an Optimal Control Problem (OCP).

The event of $x(t)$ reaching or leaving some boundary $\partial R_i$ is called a switch. 
In this paper, we consider systems with a finite number of switches on finite time intervals, i.e., Zeno trajectories are excluded.
The ODE \eqref{eq:pss} is not defined on the region boundaries $\partial R_i$, and classical notions of solutions \cite{Cortes2008a} are not sufficient to treat the rich behavior that emerges in a PSS. 
For example, during \textit{sliding modes} $x(t)$ must evolve on $\partial R_i$~\cite{Cortes2008a,Filippov1988}.
A sufficiently regular and practical notion is given by the Filippov extension for \eqref{eq:pss}.
The special structure of the PSS allows the definition of a finite number of convex multipliers $\theta_i$ and the Filippov DI reads as~\cite{Filippov1988,Stewart1990a}:
\begin{align}\label{eq:filippov_di_with_multiplers}
	\begin{split}
	\dot{x}  \in F_{\mathrm{F}}(x,u) = \Big\{ &\sum_{i\in \mathcal{J}}
			f_i(x,u) \, \theta_i  \mid \sum_{i\in \mathcal{J}}\theta_i = 1, \theta_i \geq 0,
			\\&	 \theta_i = 0 \  \mathrm{if} \;  x \notin \overline{R}_i,
			\forall  i  \in \mathcal{J} \Big\}. 		
	\end{split}
\end{align}

Let the regions $R_i$ be defined by smooth switching functions $c_j(x), j \in \mathcal{C} \coloneqq \{1,\ldots,n_c\}$.
The definition of some set $R_i$ does not have to depend on all functions $c_j(x)$.
Therefore, with $n_c$ scalar functions we define up to $n_f \leq 2^{n_c}$ regions. 
For example:
\begin{align*}
	R_1 &= \{x \in \R^{n_x}  \mid c_1(x) >0\},\\
	R_2 &= \{x \in \R^{n_x}  \mid c_1(x) <0,c_2(x) >0 \},\\
	\vdots\\
	R_{n_f} &= \{x \in \R^{n_x}  \mid c_1(x) <0,c_2(x) <0, \ldots, c_{n_c}(x) <0 \}.
\end{align*}
Note that the boundaries of the regions $\partial R_i $ are subsets of the zero-level sets of appropriate functions $c_j(x)$.
We can compactly express the definitions of the sets $R_i$ via a matrix ${S} \in \R^{\Nsys \times n_c}$, which in our example reads as:
\begin{align}\label{eq:sparse_sign_matrix}
	{S} = \begin{bmatrix}
		1 & 0 & \dots & 0\\
		-1 & 1 & \dots &  0\\
		\vdots &  \vdots & \cdots &  \vdots \\
		-1 & -1 & \dots &  -1\\
	\end{bmatrix}.
\end{align}
It is not allowed for ${S}$ to have a row with all zeros.
The sparsity in the matrix ${S}$ arises from the geometry of the regions $R_i$.
Furthermore, for every region $R_i$, we define an index set containing the indices of all switching functions relevant to its definition, i.e.,
\begin{align*}
	\mathcal{C}_i = \{ j \in  \mathcal{C} \mid {S}_{i,j}\neq 0\}, \textrm{ for all } i \in \mathcal{J}.
\end{align*}
Now, the matrix ${S}$ enables us to compactly express the definitions of the regions $R_i$ as:
\begin{align}\label{eq:standard_sets_step}
	{R}_i &= \{ x\in \R^{n_x} \mid {S}_{i,j} c_j(x)>0,\; j \in \mathcal{C}_i\}.
\end{align}
The next question to be answered is: Given the definitions of $R_i$ via the switching functions $c(x)$, how can we compute the Filippov multipliers $\theta$ in~\eqref{eq:filippov_di_with_multiplers}?
To derive such expressions, we make use of the set-valued step functions and the definition of the regions $R_i$ via the switching functions $c_j(x)$.
Let us first illustrate the development so far with an example.
\begin{example}\label{ex:filippv_set}
We regard three regions defined via the switching functions $c_1(x)$ and $c_2(x)$:
$R_1 = \{ x\in \R^{n_x} \mid c_1(x) > 0\},\;
R_2 = \{ x\in \R^{n_x} \mid  c_1(x) < 0, c_2(x)>0\}$ and
$R_3 = \{ x \in \R^{n_x}  \mid c_1(x) < 0, c_2(x)<0\}$, and the associated vector fields $f_i(x),i=1,2,3$.
By using \eqref{eq:standard_sets_step} these sets can be compactly defined via the matrix
	\begin{align*}
		{S} = \begin{bmatrix}
			1& 0\\
			-1& 1\\
			-1&-1
		\end{bmatrix}
	\end{align*}
Next, let $\alpha \in \Gamma(c(x)) \in \R^2$.
A selection of the Filippov set \eqref{eq:filippov_di_with_multiplers} and the associated ODE reads as:
\begin{align*}
		\dot{x} = \alpha_1 f_1(x)  \!+\! (1-\alpha_1)\alpha_2 f_2(x) \!+\! (1-\alpha_1)(1-\alpha_2) f_3(x).
\end{align*}
By inspection, we conclude that $\theta_1 = \alpha_1$, $\theta_2 = (1-\alpha_1)\alpha_2$ and $\theta_3 = (1-\alpha_1)(1-\alpha_2)$.
Since $\alpha_1, \alpha_2 \in [0,1]$ it is clear that $\theta_i \in [0,1], i=1,2,3$.
Similarly, by direct calculation, we verify that $\theta_1 + \theta_2 + \theta_3 = 1$.
Observe that the entries of ${S}_{i,j}$ determine how $\alpha_j$ enters the expression for~$\theta_i$.
For ${S}_{i,j} =1$ we have $\alpha_j$, for ${S}_{i,j} =-1$ we have $(1-\alpha_j)$ and for ${S}_{i,j} =0$, $\alpha_j$ does not appear in the expression for~$\theta_i$.
\end{example}
We generalize the patterns observed in our example and define the set
\begin{align}
	\label{eq:step_set}
		F_{\mathrm{S}}(x) &\! \coloneqq \!
		\Big\{ \sum_{i=1}^{\Nsys} \prod_{j \in \mathcal{C}_i} 	\Big( \frac{1-{S}_{i,j}}{2}\! + \! {S}_{i,j}\alpha_i \Big) f_i(x) \mid \alpha \! \in  \! \Gamma(c(x)) \Big\}.
\end{align}
Note that we have
\begin{align*}
		\frac{1-{S}_{i,j}}{2}+{S}_{i,j}\alpha_i	&= \begin{cases}
			\alpha_i,\ & \textrm{ if } {S}_{i,j} = 1, \\
			1-\alpha_i,\ & \textrm{ if }  {S}_{i,j} = -1.
		\end{cases}
\end{align*}
Similar definitions of the set $F_{\mathrm{S}}(x)$ can be found in \cite[Section 4.2]{Dieci2011} and \cite[Section 2.1]{Guglielmi2022}.
However, they are restricted to fully dense matrices $S$ and do not focus on developing high-accuracy discretization methods for such systems.
Next we show that $F_{\mathrm{S}}(x)$ is indeed the same set as $F_{\mathrm{F}}(x)$, i.e., the set in the r.h.s. of \eqref{eq:filippov_di_with_multiplers}.
\begin{lemma}[Lemma 1.5 in \cite{Dieci2011}]\label{lem:dieci}
		Let $a_1,a_2,\ldots a_m \in \R$. Consider the $2^m$ non-repeated products of the form
		$p_i = (1\pm a_1)(1\pm a_2)\cdots(1\pm a_m)$, then it holds that $\sum_{i}^{2^m} p_i = 2^m$.
	\end{lemma}
\begin{proposition}
		Let \begin{align}\label{eq:theta_via_step}
			\theta_i &= \prod_{j\in \C_i}  \Big(\frac{1-{S}_{i,j}}{2}+{S}_{i,j}\alpha_j\Big),
	\end{align}
	for all $i\in \mathcal{J} = \{1,\ldots,\Nsys\}$, then it holds that $F_{\mathrm{F}}(x)$ = 	$F_{\mathrm{S}}(x)$.
	\end{proposition}
	\textit{Proof.}
	We only need to show that $\theta_i \geq 0$ for all $i \in \J$ and $\sum_{i\in \J} \theta_i  =1$.
	It is easy to see that $\theta_i\in [0,1]$ as it consists of a product of terms that takes value in $[0,1]$.
	Next we show that $\sum_{i\in \mathcal{J}} \theta_i  =1$.
	We introduce the change of variables:
		$\frac{1+b_j}{2} = \alpha_j,\; \frac{1-b_j}{2}= 1-\alpha_j.$
	Then all $\theta_i$ are of the form
	\begin{align*}
		\theta_i = 2^{-|\C_i|} \prod_{j \in \C_i}  (1\pm b_j).
	\end{align*}
	If the matrix ${S}$ is dense we have that $\C_i = \C$ for all $i \in\mathcal{J}$ and $\Nsys = 2^{n_c}$.
	By applying Lemma \ref{lem:dieci} we conclude that $\sum_{i\in \mathcal{J}}\theta_i=1$ and the proof is complete.
	On the other hand, if the matrix ${S}$ has zero entries, we have that $\Nsys < 2^{n_c}$.
	We extend sequence $\{\theta_i\}_{i=1}^{n_f}$ to  $\{\tilde{\theta}_l\}_{l=1}^{2^{n_c}}$, where the terms $\tilde{\theta}_l$ are defined as follows.
	If $\C_i = \C$, then $\tilde{\theta}_i = \theta_i$.
	Now let $\C_i \subset \C$ and $\C \setminus \C_i = \{ k\}$, i.e., only one ${S}_{i,k}$ in ${S}_{i,\bullet}$ is zero.
	We can use the simple identity $\theta_i = \theta_i \frac{(1+b_k)}{2}+\theta_i \frac{(1-b_k)}{2}$, and let two additional $\tilde{\theta}_{l}$ be the two terms in the extended sum above.
	Applying this procedure inductively we obtain for all $i$ where $\C_i \subset \C$ terms $\tilde{\theta}_l$ of the form $\frac{(1\pm b_1)}{2}\cdots\frac{(1\pm b_{n_c})}{2}$.
	Now we can apply Lemma \ref{lem:dieci} and conclude that ${\sum_{i=1}^{n_f} \theta_i = \sum_{l=1}^{2^{n_c}} \tilde{\theta}_l = 1}$.
	\qed
	
\subsection{The equivalent dynamic complementarity system}	
Next, we pass from the abstract definition of the Filippov systems via set-valued step functions to the computationally more practical formulation of a DCS.
The complementarity conditions encode all the combinatorial structure and nonsmoothness in the system but can still be efficiently treated via derivative-based optimization methods~\cite{Scholtes2001,Anitescu2007}.

To perform this transition, we express the set-valued step function $\Gamma(c(x))$ as the solution map of a linear program parametric in $x$~\cite{Acary2010,Acary2014}:
\begin{subequations}\label{eq:step_parametric_lp}
	\begin{align}
		\Gamma(c(x))  = \arg \min_{\alpha\in \R^{n_c}} \ &-c(x)^\top \alpha \\
		\quad \mathrm{s.t.} \quad & 0 \leq \alpha_i \leq 1, \ i=1,\ldots,n_c \label{eq:step_parametric_lp_ineq}.
	\end{align}
\end{subequations}
Let $\lambdan,\lambdap \in \R^{n_c}$ be the Lagrange multipliers for the lower and upper bound on $\alpha$ in \eqref{eq:step_parametric_lp_ineq}, respectively.
The KKT conditions of \eqref{eq:step_parametric_lp} read as
\begin{subequations}\label{eq:step_parametric_lp_kkt}
	\begin{align}
		&c(x) = \lambdap - \lambdan, \label{eq:step_parametric_lp_kkt_lagrangian}\\
		&0 \leq  \lambdan \perp \alpha \geq 0, \label{eq:step_parametric_lp_cc1}\\
		&0 \leq  \lambdap \perp e-\alpha \geq 0 \label{eq:step_parametric_lp_cc2}.
	\end{align}
\end{subequations}
We look closer at a single component $\alpha_j$ and the associated function $c_j(x)$.
From the LP \eqref{eq:step_parametric_lp} and its KKT conditions, one can see that for $c_j(x) > 0$, we have $\alpha_j=1$.
From \eqref{eq:step_parametric_lp_kkt_lagrangian} and the complementarity condition \eqref{eq:step_parametric_lp_cc1} it follows that $\lambda_{\mathrm{p},j} = c_j(x) > 0$.
The lower bound is inactive, thus, $\lambdan_j = 0$. 
Similarly, for $c_j(x) < 0$, if follows that $\alpha_j=0$, $\lambdap_j = 0$ and $\lambdan_j = -c_j(x) > 0$.
Lastly, $c_j(x) = 0$ implies that $\alpha_j \in [0,1]$ and $\lambdap_j = \lambdan_j = 0$.
From these discussions, it is clear that $c(x)$, $\lambdan$ and $\lambdap$ are related by following expressions:
\begin{align}\label{eq:step_continuity_lambda}
	\lambdap &= \max(c(x),0),\	\lambdan = -\min(c(x),0).
\end{align}
That is, $\lambdap$ collects the positive parts of $c(x)$ and $\lambdan$ the absolute value of the negative parts of $c(x)$.
From the continuity of $c(x(t))$, it follows that the functions $\lambdap(t)$ and $\lambdan(t)$ are continuous in $t$ as well.

Using KKT systems \eqref{eq:step_parametric_lp_kkt} and combining this with the definition of the Filippov set in \eqref{eq:step_set} and the expression for $\theta_i$ in \eqref{eq:theta_via_step}, we obtain the following DCS:
		\begin{subequations}\label{eq:step_dcs}
			\begin{align}
				&\dot{x} = F(x,u)\; \theta, \\
				&\theta_i = \prod_{j\in \C_i}  \frac{1-{S}_{i,j}}{2}+{S}_{i,j}\alpha_j,\; \textrm{for all } i\in \J,\\
				&c(x) = \lambdap - \lambdan,\label{eq:step_dcs_lp_lagrangian}\\
				&0 \leq  \lambdan \perp \alpha \geq 0, \label{eq:step_dcs_cc1}\\
				&0 \leq  \lambdap \perp e-\alpha \geq 0, \label{eq:step_dcs_cc2}
			\end{align}
		\end{subequations}
		where $F(x) = [f_1(x,u),\dots, f_{\Nsys}(x,u)] \in \R^{n_x \times \Nsys}$, $\theta =(\theta_1,\dots,\theta_{\Nsys}) \in \R^{\Nsys}$ and $\lambdap,\lambdan,\alpha \in \R^{n_c}$.
		We group all algebraic equations into a single function and use a C-function $\Psi(\cdot,\cdot)$ for the complementarity condition to obtain a more compact expression:
		\begin{align*}
			G(x,\theta,\alpha,\lambdap,\lambdan) \coloneqq
			\begin{bmatrix}
				\theta_1 - \prod_{j \in C_1}  \frac{1-{S}_{1,j}}{2}+{S}_{1,j}\alpha_j\\
				\vdots\\
				\theta_{n_f} - \prod_{j \in C_{n_f}}  \frac{1-{S}_{n_f,j}}{2}+{S}_{n_f,j}\alpha_j\\
				c(x) - \lambdap + \lambdan\\
				\Psi( \lambdan,\alpha)\\
				\Psi( \lambdap,e-\alpha)
			\end{bmatrix}.
		\end{align*}
		Finally, we obtain a compact representation of \eqref{eq:step_dcs} in the form of a nonsmooth DAE:
		\begin{subequations}\label{eq:dcs_step_2}
			\begin{align}
				\dot{x} & = F(x,u)\theta,\\
				0&= 	G(x,\theta,\alpha,\lambdap,\lambdan) \label{eq:dcs_step_2_alg}.
			\end{align}
		\end{subequations}
In Table~\ref{tab:step_summary} we summarize the elementary algebraic expressions for the multipliers $\theta_i$ depending on the geometric definition of the regions $R_i$. 
Thereby, we regard the two sets $A = \{x \in \R^{n_x} \mid c_1(x) >0\}$ and $B = \{x \in \R^{n_x} \mid c_2(x) >0\}$.
All other more complicated expressions can be obtained by combining these elementary operations.
\section{Finite Elements with Switch Detection}\label{sec:fesd_step}
\subsection{Standard Runge-Kutta discretization}\label{sec:step_dcs_rk}
As a starting point for our analysis, we regard a standard RK discretization for the nonsmooth DAE formulation of the DCS \eqref{eq:dcs_step_2}.
We remind the reader that \eqref{eq:dcs_step_2_alg} collects all algebraic equations including the complementarity conditions \eqref{eq:step_dcs_cc1}-\eqref{eq:step_dcs_cc2}.
\begin{table}[t]
	\centering
	\caption{Expressions of $\theta_i$ for different definitions of $R_i$.}
	\begin{tabular}{cc}
		\hline
		\text{Definition} $R_i$ & \text{Expression} $\theta_i$\\
		\hline
		$R_i = A$ & $\theta_i=\alpha_1$ \\
		$R_i = A\cup B $& $\theta_i=\alpha_1 + \alpha_2$\\
		$R_i = A\cap B$ & $\theta_i=\alpha_1  \alpha_2$\\
		$R_i = \mathrm{int}(\R^{n_x} \setminus A) = \{x \mid c_1(x)<0\}$ & $\theta_i=1-\alpha_1$ \\
		$R_i = A \setminus B$ & $\theta_i=\alpha_1 - \alpha_2$\\
		\hline
	\end{tabular}
\vspace{-0.4cm}
	\label{tab:step_summary}
\end{table}
For ease of exposition, we regard a single control interval $[0,T]$ with a fixed control input $q \in \R^{n_u}$, i.e., we set $u(t) = q$ for $ t\in [0,T]$.
In Section \ref{sec:fesd_ocp}, we will treat the discretization of OCPs with multiple control intervals.
Let $x(0) = s_0$ be the initial value.
The control interval $[0,T]$ is divided into into $\NFE$ {finite elements} (integration intervals) $[t_{n},t_{n+1}]$ via the grid points $0= t_0 < t_1 < \ldots <t_{\NFE} = T$.
In each finite elements we regard an $\Nstg$-stage RK method which is characterized by the Butcher tableau entries $a_{i,j} ,b_i$ and $c_i$ with $i,j\in\{1,\ldots,\Nstg\}$~\cite{Hairer1991}.
The step sizes are denoted by $h_{n} = t_{n+1} - t_{n},\; n = 0, \ldots,\NFE-1$.
The approximation of the differential state at the grid points $t_n$ is denoted by $x_n \approx x(t_n)$.

We regard the differential representation of the RK method.
Hence, the derivatives of states at the stage points $t_{n,i} \coloneqq t_n + c_i h_n,\; i = 1,\ldots, \Nstg$, are degrees of freedom.
For a single finite element, we group them in the vector $V_n \coloneqq (v_{n,1}, \ldots, v_{n,\Nstg}) \in \R^{\Nstg n_x}$.
Similarly, the stage values for the algebraic variables are collected in the vectors:
$\Theta_n \coloneqq (\theta_{n,1}, \ldots, \theta_{n,\Nstg} )\in \R^{\Nstg \cdot \Nsys}$,
${A}_n \coloneqq (\alpha_{n,1}, \ldots, \alpha_{n,\Nstg} )\in \R^{\Nstg \cdot n_c}$,
$\Lambdap_n \coloneqq (\lambdap_{n,1}, \ldots, \lambdap_{n,\Nstg} )\in \R^{\Nstg \cdot n_c}$ and
$\Lambdan_n \coloneqq (\lambdan_{n,1}, \ldots, \lambdan_{n,\Nstg} )\in \R^{\Nstg \cdot n_c}$.
We collect all \textit{internal} variables in the vector
$Z_n =(x_n,\Theta_n,A_n,\Lambdap_n,\Lambdan_n,V_n)$.

The vector $x_n^{\mathrm{next}}$ denotes the value at $t_{n+1}$, which is obtained after a single integration step.
Now, we can state the RK equations for the DCS \eqref{eq:dcs_step_2} for a single finite element as
\begin{align}\label{eq:dcs_irk_single_step_reformulation}
	&0 = G_{\irk}(x_n^{\mathrm{next}}\!,Z_n,h_n,q)\!\coloneqq\! \\ \nonumber
	&
	\begin{bmatrix}
		\! v_{n,1}\! -\!  F(x_n +h_n \sum_{j=1}^{\Nstg} a_{1,j} v_{n,j},q)\theta_{n,1}\\
		\vdots\\
		v_{n,\Nstg} \! -\! F(x_n +h_n \sum_{j=1}^{\Nstg} a_{\Nstg,j} v_{n,j},q)\theta_{n,\Nstg}\\
		G(x_n + h_n\sum_{j=1}^{\Nstg} a_{1,j} v_{n,j},\theta_{n,1},\alpha_{n,1},\lambdap_{n,1},\lambdan_{n,1})\\
		\vdots\\
		G(x_n \!+\! h_n\sum_{j=1}^{\Nstg}\!a_{\Nstg,j} v_{n,j},\theta_{n,\Nstg},\alpha_{n,\Nstg},\lambdap_{n,\Nstg},\lambdan_{n,\Nstg})\\
		x_n^{\mathrm{next}} - x_n - h_n \sum_{i=1}^{\Nstg} b_i v_{n,i}
	\end{bmatrix}.
\end{align}
Next, we summarize the equations for all $\NFE$ finite elements over the entire interval $[0,T]$ in a discrete-time system format.
To simplify the statement, we need additional shorthand notation to collect all variables, on all finite elements, within the regarded control interval:
$\mathbf{x}= (x_0,\ldots,x_{\NFE}) \in \R^{(\NFE+1) n_x}$,
$\mathbf{V} = (V_0,\ldots,V_{\NFE-1}) \in \R^{\NFE \Nstg n_x}$ and
$\mathbf{h}\coloneqq (h_0,\ldots,h_{\NFE-1})\in \R^{\NFE}$.
Recall that the simple continuity condition $x_{n+1} = x_{n}^{\mathrm{next}}$ holds.
We collect the stage values of the Filippov multipliers in the vector
$\mathbf{\Theta} = ({\Theta}_0,\ldots,\Theta_{\NFE-1})\in \R^{n_{{\theta}}}$ and
$n_{{\theta}}= \NFE\Nstg\Nsys$.
Similarly, we collect the stage values of the algebraic variables specific to the step representation in vectors
$\mathbf{A},\mathbf{\Lambda}^{\mathrm{p}}, \mathbf{\Lambda}^{\mathrm{n}}\in \R^{n_{\alpha}}$, where $n_{\alpha} = \NFE\Nstg n_c$.
Finally, we collect all internal variables in the vector
$\mathbf{Z} = (\textbf{x},\mathbf{V},\mathbf{\Theta},\mathbf{A},\mathbf{\Lambda}^{\mathrm{p}},\mathbf{\Lambda}^{\mathrm{n}})\in \R^{n_{\mathbf{Z}}}$, where $n_{\mathbf{Z}} = (\NFE+1)n_x + \NFE\Nstg n_x + n_{\theta}+3n_{\alpha}$.

All computations over a single control interval of the {standard RK discretization} are summarized in:
\begin{subequations}\label{eq:dcs_irk_step_reformulation}
	\begin{align}
		{s}_1 \! &= \!F_{\mathrm{std}}(\textbf{Z}),\\
		0\! &=\! {G}_{\mathrm{std}}(\mathbf{Z},\mathbf{h},s_0,q),
	\end{align}
\end{subequations}
where $s_1\in \R^{n_x}$ is the approximation of $x(T)$ and
\begin{align*}
	F_{\mathrm{std}}(\textbf{Z})  &= x_{\NFE},\\
	G_{\mathrm{std}}(\mathbf{Z},\mathbf{h},s_0,q)
	\coloneqq	&
	\begin{bmatrix}
		x_0- s_0\\
		G_{\irk}(x_1,Z_0,h_0,q)\\
		\vdots\\
		G_{\irk}(x_{\NFE},Z_{\NFE-1},h_{\NFE-1},q)
	\end{bmatrix}.
\end{align*}
In \eqref{eq:dcs_irk_step_reformulation}, $\mathbf{h}$ is a given parameter and implicitly fixes the discretization grid.
We proceed by letting $\mathbf{h}$ be degrees of freedom and introduce the cross complementarity conditions.
\subsection{Cross complementarity}\label{sec:fesd_cross_comp_step_reformulation}
For brevity, we regard in this paper only RK methods with $c_{\Nstg} =1$, which already covers many schemes, e.g., Radau IIA and several Lobatto methods~\cite{Hairer1991}.
This means that the right boundary point of a finite element is a stage point since $t_{n+1} = t_n+c_{\Nstg} h_n$.
We will provide extensions for $c_{\Nstg} \neq 1$ in future work.

As in any event based method, we assume that there is a finite number of switches. 
To be able to detect all switches, we assume that $\NFE$ is greater then total the number of switches.
Our goal is to derive additional constraints that will allow active-set changes only at the boundary of a finite element.
Moreover, in this case, the step size $h_n$ should adapt such that all switches are detected exactly.
Note that in the standard discretization, at every RK-stage point, we have for $n=1,\ldots,\NFE$, the complementarity conditions:
\begin{subequations}\label{eq:fesd_standard_comp_step_reformulation}
	\begin{align}
		0 \leq &\lambdan_{n,m} \perp \alpha_{n,m} \geq  0,\; m =1,\ldots,\Nstg,\\
		0 \leq &\lambdap_{n,m} \perp e-\alpha_{n,m} \geq  0,\; \
		\;  m =1,\ldots,\Nstg. 
	\end{align}
\end{subequations}
As a first step, we exploit the continuity of the Lagrange multipliers $\lambdap$ and $\lambdan$.
For this purpose, we regard the boundary values of the approximation of $\lambdap$ and $\lambdan$ on an interval $[t_n,t_{n+1}]$, which are denoted by
$\lambdap_{n,0},\; \lambdan_{n,0}$ at $t_n$ and $\lambdap_{n,\Nstg},\; \lambdan_{n,\Nstg}$ at $t_{n+1}$.
We impose a continuity condition for the discrete-time versions of $\lambdap$ and $\lambdan$
for $n= 0,\ldots, \NFE-1$:
\begin{align}\label{eq:continuity_of_lambda_step_representation}
	\lambdap_{n,\Nstg}= \lambdap_{n+1,0},\;
	\lambdan_{n,\Nstg}= \lambdan_{n+1,0},\;
\end{align}
In the sequel, we use only the right boundary points $\lambdap_{n,\Nstg}$ and $\lambdan_{n,\Nstg}$, which are for $c_{\Nstg}=1$, already variables in the RK equations \eqref{eq:dcs_irk_step_reformulation}.
\begin{remark}\label{rem:lambda00_step_reformulation}
	It is important to note that $\lambdap_{0,0}$ and $\lambdan_{0,0}$ are not defined via Eq.~\eqref{eq:continuity_of_lambda_step_representation}, as we do not have a preceding finite element for $n=0$.
	However, they are crucial for determining the active set in the first finite element.
	They are not degrees of freedom but can be pre-computed for a given $x_0$.
	Using equation \eqref{eq:step_continuity_lambda} we have $\lambdap_{0,0} = \max(c(x_0),0)$ and ${\lambdan_{0,0}= -\min(c(x_0),0)}$.
\end{remark}
At a switch of the PSS, i.e., at an active-set change in the DCS \eqref{eq:step_dcs}, we have $c_i(x) = 0$.
From Eq. \eqref{eq:continuity_of_lambda_step_representation} and, due to continuity, it follows that $\lambdap_i(t)$ and $\lambdan_i(t)$ must be zero at an active-set change, as well.
Moreover, on an interval $t\in(t_n,t_{n+1})$ with a fixed active set, the components of these multipliers are either zero or positive on the whole interval.
We must now impose that their discrete-time counterparts, i.e., the stage values $\lambdap_{n,m}$ and $\lambdan_{n,m}$, have similar properties.
We achieve this with the cross complementarity conditions, which read for  $n = 0,\ldots,\NFE\!-\!1$, $m =1,\ldots, \Nstg,\  m'=0,\ldots, \Nstg$, and $m  \neq m'$ as:
\begin{subequations}\label{eq:cross_cc_true_step_reformulation}
	\begin{align}
		&0  = \mathrm{diag}(\lambdan_{n,m'})\alpha_{n,m}, 
		\\
		\begin{split}
			&0  = \mathrm{diag}(\lambdap_{n,m'})(e-\alpha_{n,m}), %
		\end{split}		
	\end{align}
\end{subequations}
In contrast to Eq.~\eqref{eq:fesd_standard_comp_step_reformulation}, we have conditions relating variables corresponding to different RK stages within a finite element.

We formalize the claims about the constraints~\eqref{eq:cross_cc_true_step_reformulation} in the next lemma.
Recall that in our notation, $\alpha_{n,m,j}$ is the $j$-th component of the vector $\alpha_{n,m}$.
\begin{lemma}\label{lem:cross_cc_statemnt_step_representation}
	Regard a fixed $n \in \{0,\ldots,\NFE\!-\!1\}$ and a fixed $j \in \C$.
	If any $\alpha_{n,m,j}$ with $m \in \{1,\ldots, \Nstg\}$ is positive, then all $\lambdan_{n,m',j}$ with $m'\in \{0,\ldots, \Nstg\}$ must be zero.
	Conversely, if any $\lambdan_{n,m',j}$ is positive, then all $\alpha_{n,m,j}$ are zero.
\end{lemma}
\textit{Proof.} Let $\alpha_{n,m,i}$ be positive, and suppose $\lambdan_{n,j,i} = 0 $ and $\lambdan_{n,k,i}>0$ for some $k,j\in \{0,\ldots, \Nstg\}, k\neq j$, then $\alpha_{n,m,i}\lambdan_{n,k,i} >0$ which violates \eqref{eq:cross_cc_true_step_reformulation}, thus all $\lambdan_{n,m',i}=0,\ m'\in \{0,\ldots, \Nstg\}$.
The converse is proven similarly. \qed

An analogous statement holds for $\lambdap_{n,m}$ and $(e-\alpha_{n,m})$.
For the switch detection, it is crucial to include the boundary points of the previous finite element in the cross complementarity conditions~\eqref{eq:cross_cc_true_step_reformulation}, namely $\lambdap_{n+1,0} = \lambdap_{n,0}$ and $\lambdan_{n+1,0} = \lambdan_{n,0}$.
A consequence of Lemma \ref{lem:cross_cc_statemnt_step_representation} is that, if the active-set changes in the $j$-th component between the $n$-th and $n+1$-st finite element, then it must hold that $\lambdap_{n,\Nstg,j} =   \lambdap_{n+1,0,j} = 0$ and $\lambdan_{n,\Nstg,j} =  \lambdan_{n+1,0,j} = 0$.
Since $x_n^{\mathrm{next}} = x_{n+1}$, we have from \eqref{eq:step_parametric_lp_kkt_lagrangian} and \eqref{eq:dcs_irk_single_step_reformulation} the condition
\begin{align*}
	c_j(x_{n+1}) = 0,
\end{align*}
which defines the switching surface between two regions.
Therefore, we have implicitly a constraint that forces $h_n$ to adapt such that the switch is detected exactly.

For clarity, the conditions \eqref{eq:cross_cc_true_step_reformulation} are given in their sparsest form.
However, the nonnegativity of $\alpha_{n,m},\lambdap_{n,m}$ and $\lambdan_{n,m}$ allows many equivalent and more compact forms.
For instance, we can use inner products instead of component-wise products, or we can even summarize all constraints for a finite element or all finite elements in a single equation, cf. \cite{Nurkanovic2022,Nurkanovic2022b} for a similar discussion.
We collect the conditions \eqref{eq:cross_cc_true_step_reformulation} into the equation $G_{\mathrm{cross}}(\mathbf{A},\mathbf{\Lambda^{\mathrm{p}}},\mathbf{\Lambda^{\mathrm{n}}}) = 0$. 
\subsection{Step equilibration}\label{sec:fesd_step_equilibration_step_reformulation}
To complete the derivation of the FESD method for \eqref{eq:step_dcs}, we need to derive the step equilibration conditions.
If no active-set changes happen, the cross complementarity constraints \eqref{eq:cross_cc_true_step_reformulation} are implied by the standard complementarity conditions~\eqref{eq:fesd_standard_comp_step_reformulation}.
Therefore, we end up with a system of equations with more degrees of freedom than conditions.
The step equilibration constraints aim to remove the degrees of freedom in the appropriate $h_n$ if no switches happen.
We achieve the goals outlined above via the equation:
\begin{align}\label{eq:step_eq_step_reformulation}
	\begin{split}
	0&= G_{\mathrm{eq}}(\mathbf{h},\mathbf{A},\mathbf{\Lambda^{\mathrm{p}}},\mathbf{\Lambda^{\mathrm{n}}})
	\coloneqq 
	\\
	&\begin{bmatrix}
		(h_{1}-h_{0})\eta_1(\mathbf{A},\mathbf{\Lambda^{\mathrm{p}}},\mathbf{\Lambda^{\mathrm{n}}}) \\
		\vdots\\
		(h_{\NFE\!-\!1}-h_{\NFE\!-\!2})\eta_{\NFE\!-\!1}(\mathbf{A},\mathbf{\Lambda^{\mathrm{p}}},\mathbf{\Lambda^{\mathrm{n}}})
	\end{bmatrix},
	\end{split}
\end{align}
where $\eta_{n}$ is an indicator function that is zero only if a switch occurs, otherwise its value is strictly positive.
In other words, if a switch happens, the $n$-th condition in \eqref{eq:step_eq_step_reformulation} is trivially satisfied.
Otherwise, it provides a condition that removes the spurious degrees of freedom.
For brevity, we omit to derive the expressions for $\eta_{n}$.
They can be obtained by similar reasoning as in~\cite[Section 3.2.3]{Nurkanovic2022}.
\subsection{Summary of the FESD discretization}\label{sec:fesd_all_fesd_eq_step}
We have now introduced all extensions needed to pass from a standard RK~\eqref{eq:dcs_irk_step_reformulation} to the FESD discretization.
With a slight abuse of notation, we collect all equations in a discrete-time system form:
\begin{subequations} \label{eq:fesd_compact_step_representation}
	\begin{align}
		s_{1} \!&=\! F_{\fesd}(\mathbf{Z}), \label{eq:fesd_compact_state_transition_step_representation}\\
		0 \!&= \!G_{\fesd}(\mathbf{Z},\mathbf{h},s_0, q , T),
	\end{align}
\end{subequations}
where $F_{\fesd}(\mathbf{x})\!=x_{\NFE}$ is the state transition map and  $G_{\fesd}(\mathbf{x},\mathbf{h},\mathbf{Z},q, T)$ collects all other internal computations including all RK steps within the regarded time interval:
\begin{align}\label{eq:fesd_compact_algebraic_step_representation}
	&G_{\fesd}(\mathbf{Z},\mathbf{h},s_0,q, T)\coloneqq
	\begin{bmatrix}
		{G}_{\mathrm{std}}(\mathbf{Z},\mathbf{h},s_0,q,T)\\
		G_{\mathrm{cross}}(\mathbf{A},\mathbf{\Lambda^{\mathrm{p}}},\mathbf{\Lambda^{\mathrm{n}}})\\
		G_{\mathrm{eq}}(\mathbf{h},\mathbf{A},\mathbf{\Lambda^{\mathrm{p}}},\mathbf{\Lambda^{\mathrm{n}}})\\
		\sum_{n=0}^{\NFE-1} h_n - T
	\end{bmatrix}.
\end{align}
Here, the control variable $q$, horizon length $T$, and initial value $s_0$ are given parameters, but $\mathbf{h}$ are degrees of freedom.

\subsection{Direct optimal control with FESD}\label{sec:fesd_ocp}
Next, we discretize this OCP using the FESD method.
The discretization process is fully automated within \texttt{NOSNOC}~\cite{Nurkanovic2022b}.
Consider $\Nctrl\geq 1$ control intervals of equal length, indexed by $k$.
We take piecewise constant control discretization, where the control variables are collected $\mathbf{q} = (q_0,\ldots,q_{\Nctrl-1})\in \R^{\Nctrl n_u}$.
All considerations can be easily extended to different control parametrizations.
We add the index $k$ to all internal variables.
On each control interval $k$, we use the FESD discretization \eqref{eq:fesd_compact_step_representation} with $N_{\mathrm{FE}}$ internal finite elements.
The state values at the control interval boundaries are grouped in the vector  $\mathbf{s} = (s_0,\ldots,s_{\Nctrl})\in\R^{(\Nctrl+1)n_x}$.
In ${\mathcal{Z}} = (\mathbf{{Z}}_0,\ldots,\mathbf{{Z}}_{\Nctrl-1})$ all internal variables, and in $\mathcal{H} = (\mathbf{h}_0,\ldots,\mathbf{h}_{\Nctrl-1})$ we collect all step sizes.

The discrete-time variant of \eqref{eq:ocp} read as:
\begin{subequations}\label{eq:ocp_discrete_time}
	\begin{align}
		\min_{\mathbf{s},\mathbf{q},\mathcal{Z},\mathcal{H}} \quad & \sum_{k=0}^{\Nctrl-1} \hat{L}(s_k,\mathbf{x}_k,q_k)+ {R}(s_{\Nctrl}) \\
		\textrm{s.t.} \quad  &s_{0} = \bar{x}_0,\\
		&{s}_{k+1}  = F_{\fesd}(\mathbf{x}_k),\;  k = 0,\ldots,\Nctrl\!-\!1,\\
		&0 = G_{\fesd}(\mathbf{x}_k,\mathbf{Z}_k,q_k),\; \! k = 0,\ldots,\Nctrl\!\!-\!\!1,\\
		&0 \geq G_{\mathrm{p}}(s_k,q_k),\; k = 0,\ldots,\Nctrl-1,\\
		&0 \geq G_{\mathrm{t}}(s_{\Nctrl}),
	\end{align}
\end{subequations}
where $\hat{L}:\R^{n_x}\times \R^{(\NFE+1)\Nstg n_x} \times \R^{n_u}\to \R$ is the discretized running costs.
Due to the complementarity constraints in the FESD discretization, \eqref{eq:ocp_discrete_time} is an MPCC.
In practice, MPCCs can usually be solved efficiently by solving a sequence of related and relaxed NLPs within a homotopy approach.
Such an approach, with some of the standard reformulations~\cite{Scholtes2001,Anitescu2007}, is implemented in \texttt{NOSNOC}.
The underlying NLPs are solved via \texttt{IPOPT}~\cite{Waechter2006} called via its \texttt{CasADi} interface~\cite{Andersson2019}.
\begin{figure}
	\includegraphics[width=0.90\columnwidth]{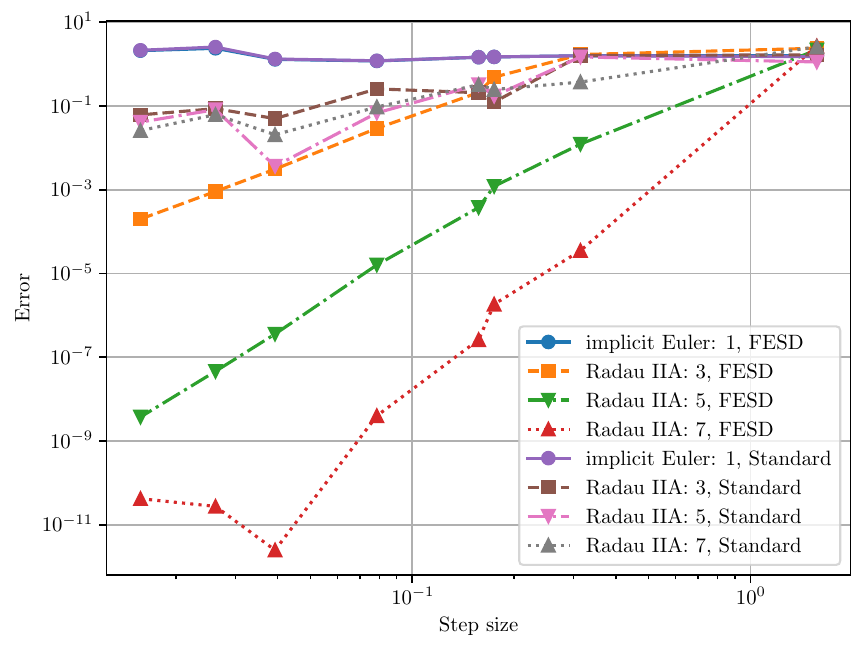}
	\caption{Accuracy vs. step size: Simulation of example \eqref{eq:ode_oscilator} with different RK schemes and step sizes.
		The number next to the method's name is the order of the underlying RK method.}
	\label{fig:fesd_integration_order}
	\vspace{-0.55cm}
\end{figure}
\section{Numerical examples}
We show on a numerical simulation example that FESD recovers the order of accuracy that RK methods have for smooth ODEs. 
We further use a hopping robot OCP example modeled as a system with state jumps in order to show that the step reformulation can outperform Stewart's reformulation used in~\cite{Nurkanovic2022}. 
\begin{figure*}[t]
	\centering 
		\includegraphics[width=0.9\linewidth]{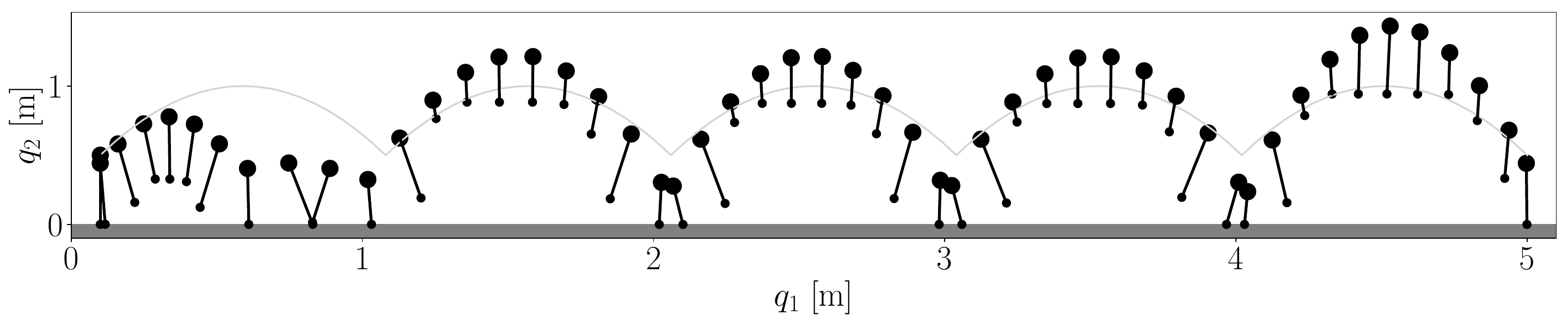}
	\caption{Several frames of the hopping robot trajectory.}
	\label{fig:solution_ocp}
		\vspace{-0.5cm}
\end{figure*}
\begin{figure}[t]
	\centering 
	\includegraphics[width=0.95\linewidth]{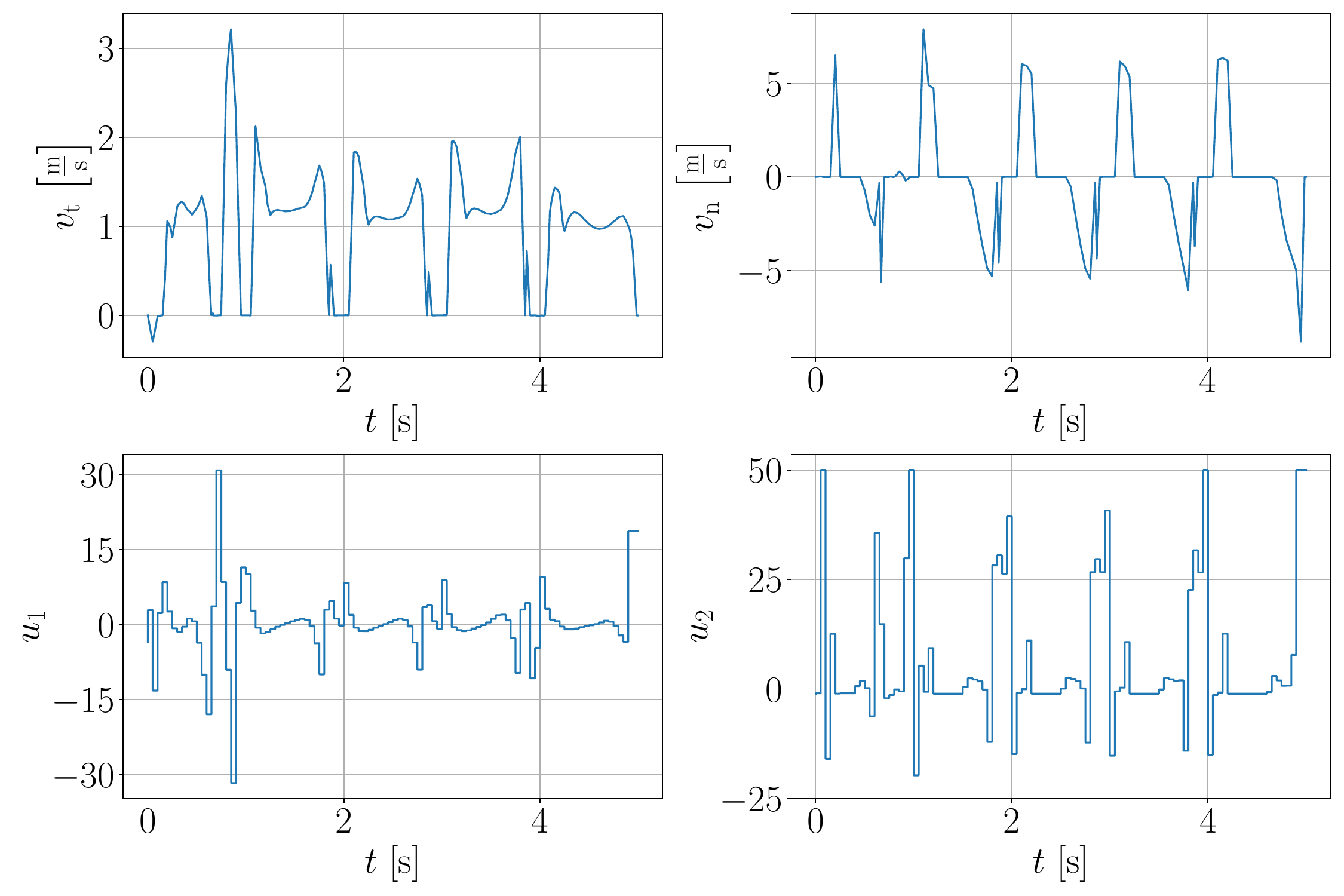}
	\caption{The top plots show the discontinuous normal and tangential velocity plotted over the \textit{physical time} $t$, the bottom plots the optimal controls.}
	\label{fig:velocities_controls}
	\vspace{-0.65cm}
\end{figure}

\subsection{Integration order experiment}
We compare the standard RK time-stepping method for DCS from Eq. \eqref{eq:dcs_irk_step_reformulation} to FESD \eqref{eq:fesd_compact_step_representation} on a simulation example from~\cite{Nurkanovic2022}.
Regard the PSS:
\begin{align}\label{eq:ode_oscilator}
	\dot{x} = f_i(x),\ \mathrm{if}\ x \in R_i,
\end{align}
with $f_1(x) = A_1x, f_2(x) = A_2x, c(x) = \|x\|_2^2 -1$, $R_1 = \{x \mid c(x)<0 \}$, $R_2 = \{x \mid c(x)>0 \}$.
The system matrices are
%
\begin{align*}
A_1 = \begin{bmatrix}
	1 & 2 \pi \\ -2 \pi &1
\end{bmatrix},
\ A_2 = \begin{bmatrix}
	1 & -2 \pi \\ 2 \pi &1
\end{bmatrix}.
\end{align*}
The initial value is $x(0) = (e^{-1},0)$ and we regard an the interval $t \in [0,T]$ with $T = \frac{\pi}{2}$.
It can be shown that the switch happens at $\ts = 1$ and that
$x(T) = (\exp{(T \! -\! \ts)}\cos(\omega(T\! -\! \ts)),-\exp{(T \! -\! \ts)}\sin(\omega(T \! -\! \ts)))$, for $T>\ts$.
Thus, given a numerical approximation $\hat{x}(t)$, we can determine the global integration error $E(T) = \| x(T)-\hat{x}(T)\|$ and observe the accuracy order of the integrator for a varying step size $h$.

Figure~\ref{fig:fesd_integration_order} shows the results of our experiment for the Radau IIA methods.
In all cases, the standard RK time-stepping \eqref{eq:dcs_irk_step_reformulation} has, as expected, only first-order accuracy, i.e., $O(h)$.
The FESD method \eqref{eq:fesd_compact_step_representation} detects the switch and recovers the high integration order properties of the underlying RK method, i.e., $O(h^p)$ with $p=2\Nstg-1$.

\subsection{Optimal control example with state jumps}
We regard the example of a planar hopper with state jumps and friction. 
Using the time-freezing reformulation, we can transform the system with state jumps into a PSS of the form of \eqref{eq:pss}~\cite{Nurkanovic2021}.
The system is inspired by a reaction wheel and force-controlled single leg as described by~\cite{Raibert1989}. 
It is assumed that the inertia matrix is a constant diagonal matrix $M = \mathrm{diag}(m_{\mathrm{b}}+ m_{\mathrm{l}}, m_{\mathrm{b}}+ m_{\mathrm{l}},I_{\mathrm{b}}+ I_{\mathrm{l}}, m_{\mathrm{l}})$, where 
$m_{\mathrm{b}} = 1$  is the mass of the body, 
$m_{\mathrm{l}} = 0.25$ is the mass of the link, 
$I_{\mathrm{b}} = 0.25$ the inertia of the body, and
$I_{\mathrm{l}} = 0.025$ the inertia of the link.

The configuration $q=(q_1,q_2,\psi,l)$ consists of the 2D position, orientation, and leg length, respectively.
The resulting PSS state consists of the position $q(\tau)$, velocity $v(\tau)$ and clock state $t(\tau)$ (needed for time-freezing, cf. \cite{Nurkanovic2021}), i.e., $x=(q,v,t) \in \R^9$. 
The robot makes contact with the ground with the tip of its leg, and the gap function is given by 
$f_c(q) = q_2 - l\cos(\psi)$. 
We also make use of the contact normal $J_\mathrm{n}(q) = (0,1,l\sin(\psi),-\cos(\psi))$, and tangent $J_\mathrm{t}(q) = (0,1,l\cos(\psi),\sin(\psi))$.
For the time-freezing reformulation, we regard three switching functions: the gap function, as well as the normal and tangential contact velocities:
\begin{align*}
c(x) = (f_c(q), J_\mathrm{n}(q)^\top v, J_\mathrm{t}(q)^\top v).
\end{align*}
They define the three regions: 
\begin{align*}
\begin{split}
		R_1 &= \{x \in \R^{n_x} \mid f_c(q)>0\}\\
		&\cup \{x \in \R^{n_x} \mid f_c(q)<0,J_\mathrm{n}(q)^\top v>0\},
\end{split}\\
	R_2 &= \{x \in \R^{n_x} \mid f_c(q)<0,J_\mathrm{n}(q)^\top v<0,J_\mathrm{t}(q)^\top v>0\},\\
	R_3 &= \{x \in \R^{n_x} \mid f_c(q)<0,J_\mathrm{n}(q)^\top v<0,J_\mathrm{t}(q)^\top v<0\}.
\end{align*}
In region $R_1$, we define the unconstrained (free flight)  dynamics of the robot, and in regions $R_2$ and $R_3$, auxiliary dynamics that mimic state jumps in normal and tangential directions due to frictional impacts, cf.~\cite{Nurkanovic2021}:
\begin{align*}
	&f_1(x,u)  = (q, M^{-1}f_v(q,u) , 1),\\
	&f_2(x)  = (\mathbf{0}_{4,1},  M^{-1}(J_{\mathrm{n}}(q) - J_{\mathrm{t}}(q)\mu)a_{\mathrm{n}}, 0), \\
	&f_3(x)  = (\mathbf{0}_{4,1}, M^{-1}(J_{\mathrm{n}}(q) + J_{\mathrm{t}}(q)\mu)a_{\mathrm{n}} , 0),
\end{align*}
with $f_v(q,u)=(-\sin(\psi)u_2,(m_{\mathrm{b}}\!+\! m_{\mathrm{l}})g+\cos(\psi)u_2, u_1, u_2)$ summarizing all forces acting on the robot, $u =(u_1,u_2) \in \R^2$ are the controls, $\mu = 0.45$ is the coefficient of friction, $g=9.81$ the gravitational acceleration constant, and $a_{\mathrm{n}}=100$ is the auxiliary dynamics constant~\cite{Nurkanovic2021}.
Note that the clock state dynamics are $\frac{\dd t}{\dd \tau} = 1$ in $R_1$, and $\frac{\dd t}{\dd \tau} = 0$ in $R_2$ and $R_3$.
Solution trajectories of the PSS are continuous in time. However, by taking the pieces of the trajectory where $\frac{\dd t}{\dd \tau}>0$, we recover the solution of the original system, cf.~\cite{Nurkanovic2021}.
To demonstrate the efficiency gained via the step function approach, we run an experiment in which the hopper attempts to cross 5 meters with a given reference trajectory of 5 jumps in $T= 5$ seconds.
The initial value is $x(0)  = (0.1,0.5,0,0.5,0,0,0,0)$.
Given a reference $x^{\mathrm{ref}}(t)$, we define the least-squares objective with the running and terminal costs:
\begin{align*}
	L(x(\tau),u(\tau)) &=(x(\tau)-x^{\mathrm{ref}}(\tau))^\top Q (x(\tau)-x^{\mathrm{ref}}(\tau)) \\ &+ \rho_u u(\tau)^\top u(\tau),\\
	R(x(T)) &= (x(T)-x^{\mathrm{ref}}(T)^\top Q_T (x(T)-x^{\mathrm{ref}}(T)),
\end{align*}
$Q =  \mathrm{diag}(100, 100, 20, 20, 0.1, 0.1, 0.1, 0.1,0)$
$\rho_u = 0.01$, and $Q_T = \mathrm{diag}(300, 300, 300, 300, 0.1, 0.1, 0.1, 0.1,0)$.
We define the path constraints:
\begin{subequations}
\begin{align}
	& x_{\mathrm{lb}} \leq  x \leq x_{\mathrm{ub}},\\
	&u_\mathrm{lb}\leq u(t) \leq u_{\mathrm{ub}},\\
	&J_{\mathrm{t}}(q)^\top v (1-\alpha_1)(1-\alpha_2) = 0, \label{eq:robot_no_slip}
\end{align}
\end{subequations}
where $x_{\mathrm{ub}} = (5.1, 1.5, \pi, 0.5, 10, 10, 5, 5,\infty)$,
$x_{\mathrm{lb}} = (0, 0, -\pi, 0.1, -10, -10, -5, -5, -\infty)$,
$u_{\mathrm{ub}} = (50,50)$, and $u_{\mathrm{lb}} =- u_{\mathrm{ub}}$.
The last equality constraint \eqref{eq:robot_no_slip} models the following: If the normal contact force, which is proportional to $(1-\alpha_1)(1-\alpha_2)$, is nonnegative (cf. \cite{Nurkanovic2021}), then $J_{\mathrm{t}}(q)^\top v = 0$.
This prevents the optimizer from choosing controls that lead to a lot of slipping when the robot is on the ground.
Collecting all the above, we can formulate an OCP of the form of \eqref{eq:ocp}, which we discretized with the FESD Radau IIA scheme of order 3 ($\Nstg = 2$), with $\NFE =3$ finite elements on every control interval.
The OCP is discretized and solved with \texttt{nosnoc\_py} in a homotopy loop with \texttt{IPOPT}~\cite{Waechter2006}.

\begin{figure}
	\centering 
	\includegraphics[width=0.49\linewidth]{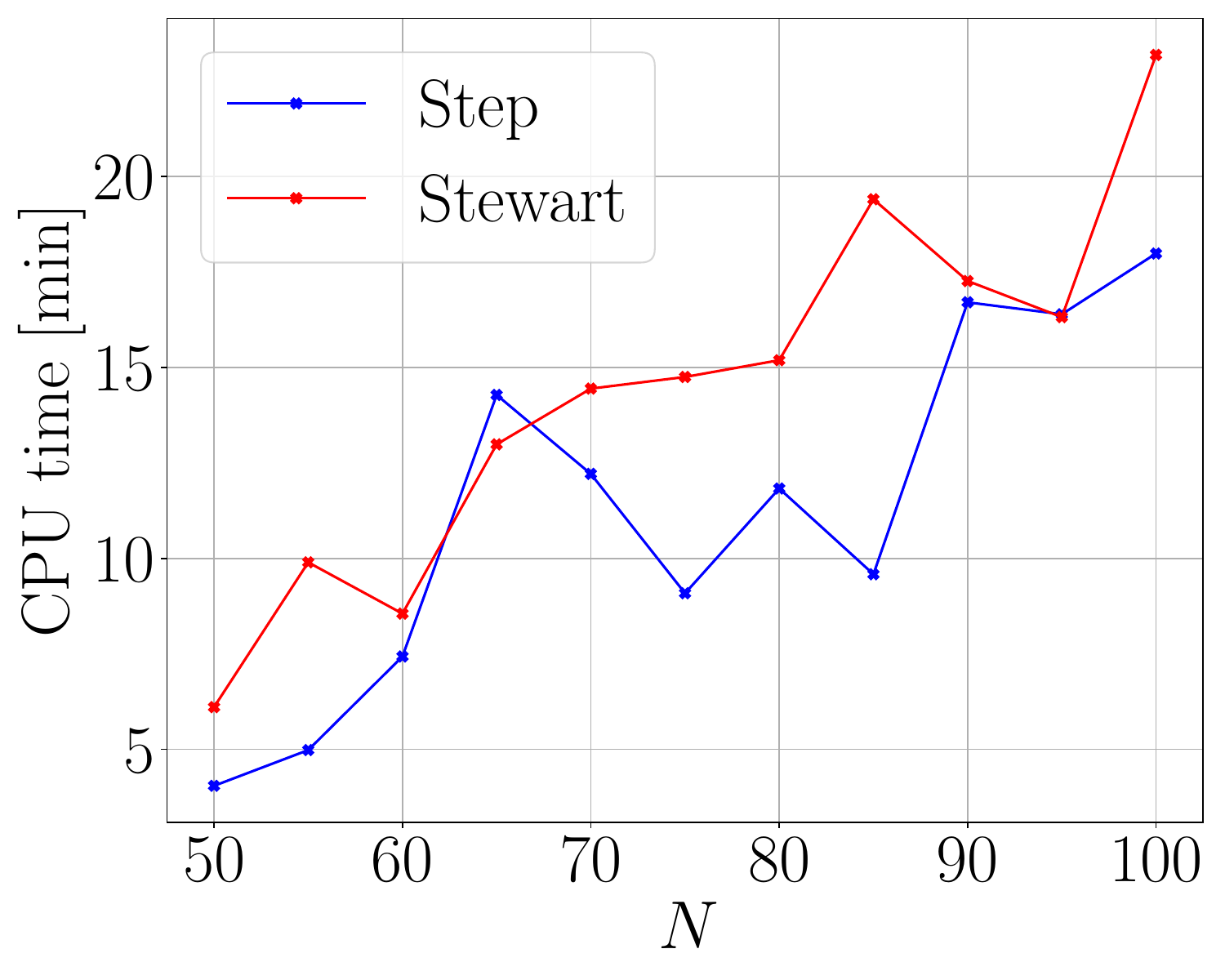}
	\vspace{-0.25cm}
	\includegraphics[width=0.49\linewidth]{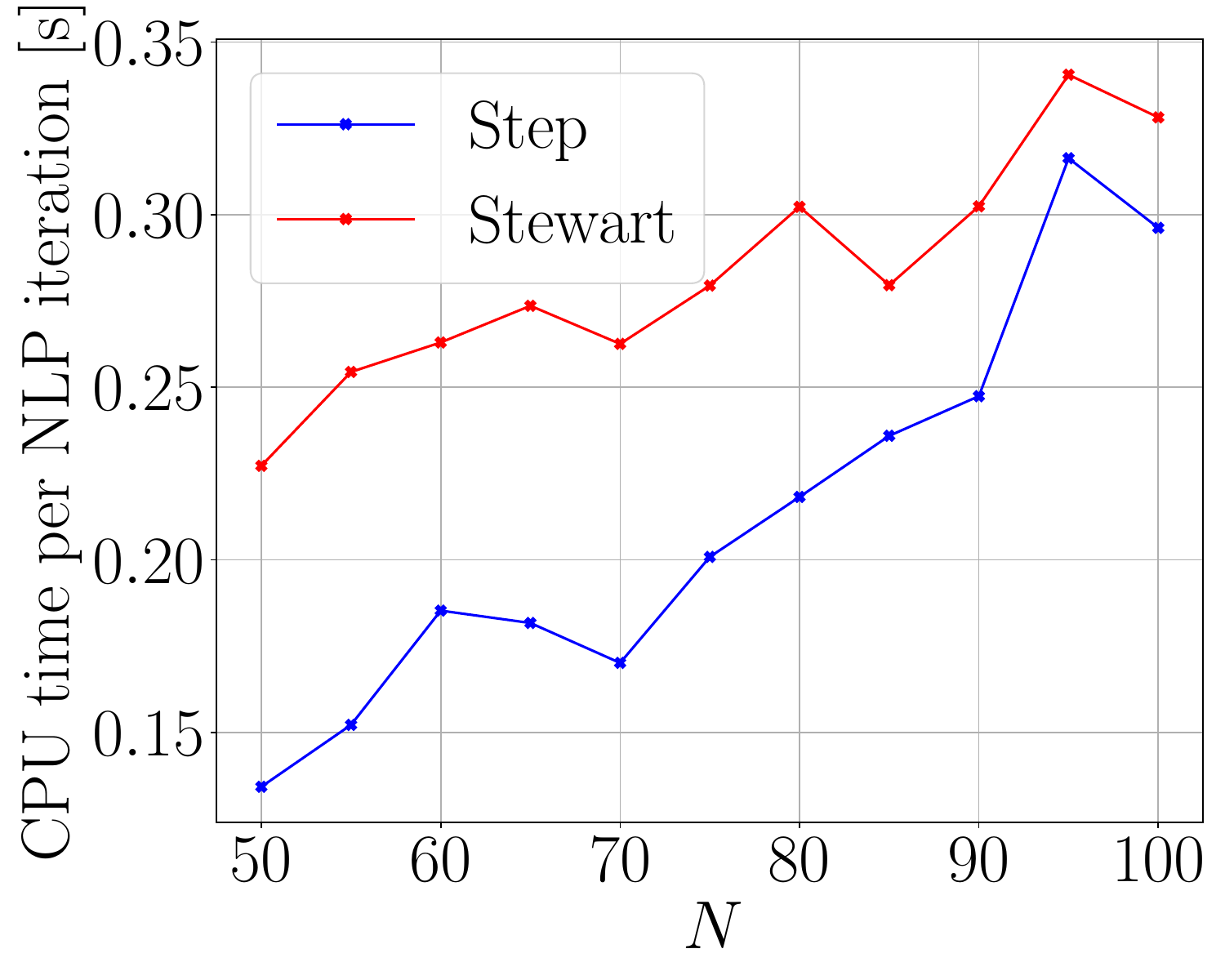}
	\caption{The total CPU time of the homotopy loop for different $N$ (left plot) and the CPU time per NLP solver iteration (right plot).}
	\label{fig:ocp_comparisson}
	\vspace{-0.55cm}
\end{figure}

Several frames of an example solution ($N=100$) can be seen in Figure~\ref{fig:solution_ocp}.
Figure~\ref{fig:velocities_controls} shows the normal and tangential velocity of the foot tip, and optimal controls. 
Next, we solve this OCP for 10 different values $N$ (number of control intervals) from 50 to 100 in increments of 5 and compare it to the FESD derived for Stewart's reformulation~\cite{Nurkanovic2022}. 
We plot the CPU time per NLP iteration and total CPU time for both approaches in Figure~\ref{fig:ocp_comparisson}. 
The step reformulation leads to faster NLP iterations than the Stewart reformulation, since it needs less variables. 
The overall computation time is governed by many factors (homotopy loop, initialization, NLP solver performance, etc.) and as such shows a less clear trend. 
However, we see that in most cases the step approach is still faster by up to 50\%. 
\section{Conclusion and outlook}
This paper introduced an extension of the Finite Elements with Switch Detection (FESD)~\cite{Nurkanovic2022} to optimal control and simulation problems with nonsmooth systems with set-valued step functions. 
This formulation covers a broad spectrum of practical problems. 
We show in numerical examples that depending on the switching functions and geometry of the underlying piecewise smooth system it can be computationally more efficient than the formulation in \cite{Nurkanovic2022}. 
Several extensions, more numerical benchmarks and a theoretical analysis of the proposed method is available in the journal version of this paper~\cite{Nurkanovic2023c}.

\bibliographystyle{plain}


\end{document}